\documentclass[10pt]{amsart}
\usepackage{latexsym, amsmath, amssymb, longtable, booktabs, amscd, graphicx}

\usepackage[english]{babel}
\usepackage[utf8]{inputenc}

\usepackage{hyperref}
\usepackage[all]{xy}
\usepackage{color}
\theoremstyle{plain}
\usepackage[margin=1in]{geometry}
\usepackage [english]{babel}
\usepackage[autostyle]{csquotes}
\usepackage{relsize}
\usepackage{dsfont}
\usepackage{bbm}
\usepackage{lipsum}
\usepackage{enumerate}   
\usepackage{enumitem}
\usepackage{tabstackengine}
\allowdisplaybreaks

\usepackage{hyperref}
\usepackage{enumerate}
\usepackage{enumitem}
\usepackage{mathtools}

\usepackage{tikz}
\usepackage{pst-node}

\usepackage[all]{xy}

\numberwithin{equation}{section}
\newcommand{\ind}{\operatorname{ind}}

\newcommand{\Q}{\mathbb{Q}}
\newcommand{\Z}{\mathbb{Z}}

\newcommand{\F}{\mathbb{F}}

\newcommand\Sym{{\mathrm {Sym}}}

\newcommand\linee{{
        {\begin{tiny}
        \begin{xymatrix}{
         \bullet  \ar@{-}[d] \\ 
         \bullet  \ar@{-}[d] \\
         \bullet} 
       \end{xymatrix}
       \end{tiny}} }}  
\newcommand\diamondd{{
        \begin{tiny}
        \begin{xymatrix}{
         & \bullet  \ar@{-}[dl] \ar@{-}[dr] &  \\ 
         \bullet  \ar@{-}[dr] &   &   \bullet  \ar@{-}[dl] \\
         & \bullet & } 
       \end{xymatrix}
       \end{tiny} }}

\newtheorem{thm}{Theorem}[section]
\newtheorem{theorem}[thm]{Theorem}
\newtheorem{cor}[thm]{Corollary}
\newtheorem{conj}[thm]{Conjecture}

\theoremstyle{definition}

\theoremstyle{remark}

\theoremstyle{remark}

\begin{document}

\title[Zig-zag and local constancy]{A zig-zag conjecture and local constancy for Galois representations}

\author[E. Ghate]{Eknath Ghate} 
\address{School of Mathematics, Tata Institute of Fundamental Research, Homi Bhabha Road, Mumbai-5, India}
\email{eghate@math.tifr.res.in}

\begin{abstract}
We make a zig-zag conjecture describing the reductions
of irreducible crystalline two-dimensional representations of $G_{\Q_p}$ of half-integral
slopes and exceptional weights. Such weights are two more than twice the
slope mod $(p-1)$. 
We show that zig-zag holds for half-integral slopes at most $\frac{3}{2}$.  We then 
explore the connection between zig-zag and local constancy results in the weight.
First we show that known cases of zig-zag force local constancy to fail for small weights. Conversely, we 
explain how local constancy forces zig-zag to fail for some small weights and 
half-integral slopes at least $2$. However, we expect zig-zag to be qualitatively true in general. 
We end with some compatibility results between zig-zag and other results.
\end{abstract}

\maketitle



\section{Introduction}

Let $p$ be an odd prime. This paper is concerned with understanding the
reductions of certain crystalline two-dimensional representations of the
local Galois group $G_{\Q_p}$. This problem is classical, and important in view of 
its applications to Galois representations attached to modular forms.  

Let $V_{k,a_p}$ be the irreducible two-dimensional crystalline representation of $G_{\Q_p}$
defined over a finite extension $E$ of $\Q_p$, of Hodge-Tate weights $(0,k-1)$ with $k \geq 2$, and positive slope $v(a_p) > 0$,  
for $a_p \in E$, and $v$  the $p$-adic valuation of $\bar{\Q}_p$ normalized so that $v(p) = 1$. The (dual of the) representation
$V_{k, a_p}$ can be described explicitly it terms of Fontaine's functor $D_\mathrm{cris}$ which sets up an equivalence
of categories between crystalline representations of $G_{\Q_p}$ over $E$ and weakly admissible filtered $\varphi$-modules
over $E$. We have
$D_\mathrm{cris}(V_{k,a_p}^*) = D_{k,a_p}$,  where
$D_{k,a_p}  = E e_1 \oplus E e_2$ 
is the filtered $\varphi$-module given by 
\begin{eqnarray*}
  \varphi(e_1) =  p^{k-1} e_2 & \text{ and } & \varphi(e_2) = - e_1 + a_p e_2, \\
  \mathrm{Fil}^i{D_{k, a_p}} & = & \begin{cases} 
                                                  D_{k,a_p}   &   \text{ if } i \leq 0, \\
                                                  E e_1         &   \text{ if } 1 \leq i \leq k-1, \\
                                                  0               &   \text{ if } i \geq k.
                                               \end{cases}
\end{eqnarray*}

Let $\bar{V}_{k,a_p}$ be the {\it semisimplification} of the 
reduction of $V_{k,a_p}$ modulo the maximal ideal of the ring of integers of $E$. It is a two-dimensional semisimple representation of
$G_{{\mathbb Q}_p}$ defined over $\bar{\mathbb F}_p$, and is independent of the choice of lattice used to define the reduction.  
If $f = \sum_{n=1}^\infty a_n q^n$ is a primitive cusp form of weight $k \geq 2$, level coprime to $p$, and (for simplicity) trivial nebentypus character,
then the local Galois representation $\rho_f |_{G_{\Q_p}}$ attached to $f$ and $p$ is isomorphic to $V_{k,a_p}$,  at least if $a_p^2 \neq 4 p^{k-1}$, 
so the reduction $\bar{\rho}_f |^{ss}_{G_{\Q_p}}$ is isomorphic to $\bar{V}_{k,a_p}$. Thus the structure of
the reductions $\bar{V}_{k,a_p}$
has important applications to the study of various objects attached to modular forms such as their  motives and
Galois representations.

It is an outstanding problem to understand the shape of the reduction $\bar{V}_{k,a_p}$.
The reduction $\bar{V}_{k,a_p}$ was computed classically by Fontaine and 
Edixhoven \cite{Edixhoven92} for all small weights $2 \leq k \leq p+1$. In a remarkable breakthrough using representation
theoretic techniques (Langlands correspondences, as developed by Breuil, Berger, Colmez, Pa\v sk\=unas, Dospinescu \cite{Bre03a}, \cite{Bre03b}, \cite{B10}, \cite{Col10}, \cite{Pas13}, \cite{CDP14}, and others), 
this range of weights was extended by Breuil \cite{Bre03b} to all $k \leq 2p+1$, at least if $p$ is odd. Using different techniques (the theory of ($\varphi, \Gamma$)-modules), Yamashita-Yasuda \cite{YY} have announced a further extension 
to weights up to $\frac{p^2+1}{2}$. 

For simplicity, let us write $v$ for the slope $v(a_p)$.
Thus, $v$ denotes both the $p$-adic valuation and the slope, depending on the context. 
The reduction $\bar{V}_{k,a_p}$ is also known for all large slopes 
$v > \lfloor \frac{k-2}{p-1} \rfloor$
by Berger-Li-Zhu \cite{BLZ}. There has been a spate of recent work  computing the reduction $\bar{V}_{k,a_p}$ for small slopes $v$. 
Buzzard-Gee \cite{BG09}, \cite{BG13} treated the case of slopes $v$ in $(0,1)$.
The case of slopes $v$ in $(1,2)$ was treated in \cite{GG15}, \cite{Bhattacharya-Ghate}, 
under an assumption when $v = \frac{3}{2}$. The case of slope $v = 1$ 
was treated in \cite{BGR18}. 

The first goal of this paper is fill a gap in the literature computing the reduction for all slopes at most $2$
by describing the complete recent treatment of the case of slope $v = \frac{3}{2}$ carried out  in \cite{GR19}. 
A more general second aim of this paper is as follows. 
Let us say that a weight $k$ is {\it exceptional} for a particular half-integral 
(and possibly integral) slope $v \in \frac{1}{2} \Z$ with $0 < v \leq \frac{p-1}{2}$ 
if $$k \equiv 2v+2 \mod (p-1).$$ 
With hindsight, it has emerged
that these weights are the hardest to treat. 
In this paper, we make a general zig-zag conjecture which describes the 
reduction $\bar{V}_{k,a_p}$ for all exceptional weights for all half-integral slopes
in a qualitative way. We also make some refinements for small half-integral slopes
that specify the reduction precisely.
Our refined conjecture specializes
to known theorems when $v = \frac{1}{2}$ \cite{BG13}, $v = 1$ \cite{BGR18}, which is not surprising since
it was modelled on these results,  
and is now known to be true for $v = \frac{3}{2}$ in view of our recent work \cite{GR19}.

As far as we are aware, the reduction problem for slope 2 is still open, though partial results for small slopes 
larger than 2 have been announced by  Arsovski 
\cite{Ars18} and Nagel-Pande \cite{NP}. In related parallel developments on 
the reduction problem, we remark that certain crystabelian cases of weight 2 and slope at most 1 had been treated earlier 
by Savitt \cite{Sav05}, and that similarly, semistable (non-crystalline) cases of small even weights $k \leq p+1$ had earlier 
been treated by Breuil-M\'ezard \cite{Breuil-Mezard02}, and more recently by Guerberoff-Park \cite{GP18} 
for the remaining odd weights in this range, using the theory of strongly divisible modules.

 
\subsection{History}
 
Before we state the zig-zag conjecture, and as motivation for it, let us describe in more detail some recent history showing how
the exceptional congruences classes of weights have emerged as the most difficult to treat. 
To this end, we recall some standard notation. Let $\omega$ and $\omega_2$
be the mod $p$ fundamental characters of levels $1$ and $2$. Let $\ind(\omega_2^c)$ be the (irreducible) mod $p$ 
representation of $G_{\Q_p}$ obtained
by inducing the $c$-th power of $\omega_2$ from the index 2 subgroup $G_{\Q_{p^2}}$ of $G_{\Q_p}$ to $G_{\Q_p}$ (for $p+1 \nmid c$). 
Let $\mu_\lambda$ be the
unramified character of $G_{{\mathbb Q}_p}$ mapping a (geometric) Frobenius at $p$ to $\lambda \in \bar{\mathbb F}_p$. 
Let $r := k-2$ and let $b \in \{1,2, \cdots, p-1\}$ represent the congruence class of $r$ 
mod $(p-1)$. Then $b = 2v$ is a representative for the exceptional congruence class of weights $r$ mod $(p-1)$. In particular,
$b = 1, 2, 3, \cdots$, represents the exceptional congruence classes of weights $r$ mod $(p-1)$ for 
the half-integral slopes $\frac{1}{2}, 1, \frac{3}{2}, \cdots$, respectively.  

In \cite{BG09}, Buzzard-Gee showed that 
the reduction $\bar{V}_{k,a_p}$ is always {\it irreducible} for slopes $v$ in $(0,1)$ (and isomorphic to  $\ind (\omega_2^{b+1})$), except possibly in the exceptional case $v = \frac{1}{2}$ and $b = 1$. This case was only treated completely in \cite{BG13}, where
the authors show that when a certain parameter, which we call $\tau$, 
is larger than another parameter, which we call $t$, a reducible possibility (namely $\omega \oplus \omega$ on inertia) occurs instead. 
More precisely, setting 
\begin{eqnarray*}
  \tau & = & v \left( \frac{a_p^2 - r p}{p a_p} \right), \\
    t  & = & v (1-r), 
\end{eqnarray*}
it is shown in \cite[Theorem A]{BG13} that in the exceptional case $v = \frac{1}{2}$ and $b = 1$, there is a {\it dichotomy}:
\begin{eqnarray*}
  \bar{V}_{k,a_p} & \sim & 
  \begin{cases}
    \mathrm{ind}(\omega_2^{b+1}),  & \text{if }  \tau < t  \\
	\mu_{\lambda} \cdot \omega^b \,\oplus\,\mu_{\lambda^{-1}} \cdot \omega, & \text{if }  \tau \geq t,
  \end{cases}  
\end{eqnarray*}
for $r > 1$, where $\lambda$ is a root of the quadratic equation  
\begin{eqnarray*}
  \lambda + \frac{1}{\lambda} & = & \overline{\frac{1}{1-r} \cdot \frac{a_p^2 - r p}{p a_p}}.
\end{eqnarray*} 

In \cite{BGR18}, the reduction $\bar{V}_{k,a_p}$ was completely determined for $p > 3$ on the boundary $v = 1$ of the Buzzard-Gee annulus $(0,1)$, 
and was shown to be generically {\it reducible} instead. In the difficult 
exceptional case $v=1$ and $b = 2$, the authors were able (after some previous iterations of the paper on the {\tt arXiv}) to
establish, for $r > 2$, a {\it trichotomy}:

%
%
%
%
\begin{eqnarray*}
  \bar{V}_{k, a_p} & \sim & 
  \begin{cases}
    \mathrm{ind}(\omega_2^{b+1}),  & \text{if } \tau < t \\
    \mu_{\lambda} \cdot \omega^b\,\oplus\,\mu_{\lambda^{-1}} \cdot \omega, & \text{if } \tau = t \\ 
    \mathrm{ind}(\omega_2^{b+p}), & \text{if } \tau > t,       
\end{cases} 
\end{eqnarray*}
where
\begin{eqnarray*}
  \tau & = & v \left( \frac{a_p^2 - \binom{r}{2} p^2}{p a_p} \right), \\
    t  & = & v (2-r), 
\end{eqnarray*}
and where $\lambda$ is a constant given by
\begin{eqnarray*}
    \lambda & = &  \overline{\dfrac{2}{2-r} \cdot \dfrac{a_p^2 -\binom{r}{2} p^2}{pa_p} }.
\end{eqnarray*}

\subsection{A conjecture} Based on these results for slopes $\frac{1}{2}$ and $1$, and some  computations of Rozensztajn \cite{Roz16} for some small half-integral
slopes, one might guess that in the general exceptional case $v \in \frac{1}{2} {\mathbb Z}$ and $b = 2v$, there are $b+1$ possibilities for
$\bar{V}_{k,a_p}$, with various {\it irreducible and reducible cases occurring alternately}.
More precisely, we make the following qualitative conjecture.

\begin{conj}[Zig-Zag Conjecture] 
  \label{zigzag}
  Say that $r = k-2  \equiv b = 2v \mod (p-1)$ is an exceptional congruence class of weights
  for a particular half-integral slope $0 < v \in \frac{1}{2} {\mathbb Z} \leq \frac{p-1}{2}$. Set $t = v(r-b)$. 
  Then, for all weights $r > b$, except possibly some $r$ that are $p$-adically close to some small weights, 
  and for all but possibly finitely many $a_p$, there is a rational     
  parameter $\tau = v(c)$, for some $c = c(r, a_p)$, such that
 as $\tau$ varies through the rational line, the (semisimplifcation of
 the) reduction $\bar{V}_{k,a_p}$ of the crystalline representation 
 $V_{k,a_p}$ satisfies a $(b+1)$-fold-{\it chotomy} consisting of alternating irreducible and reducible cases:
 \begin{eqnarray*}
  \bar{V}_{k, a_p} & \sim & 
  \begin{cases}
    \begin{array}{l}
      \mathrm{ind}(\omega_2^{b+1}),    
    \end{array}                                                              & \text{ if } \tau < t \\
    \begin{array}{l}
      \mu_{\lambda_1} \cdot \omega^b\,\oplus\,\mu_{\lambda_1^{-1}} \cdot \omega, 
    \end{array}                                                              & \text{ if } \tau = t \\ 
    \begin{array}{l}
      \mathrm{ind}(\omega_2^{b+p}),     
    \end{array}                                                              & \text{ if } t < \tau < t+1 \\
    \begin{array}{l}
      \mu_{\lambda_2} \cdot \omega^{b-1}\,\oplus\,\mu_{\lambda_2^{-1}} \cdot \omega^2, 
    \end{array}                                                              & \text{ if } \tau = t + 1 \\ 
    \begin{array}{l}
      \mathrm{ind}(\omega_2^{b+2p-1}),
    \end{array}                                                              & \text{ if } t +1  < \tau < t+2      \\
    \begin{array}{l} 
      \mu_{\lambda_3} \cdot \omega^{b-2}\,\oplus\,\mu_{\lambda_3^{-1}} \cdot \omega^3, 
    \end{array}                                                              & \text{ if } \tau = t + 2\\ 
    \>\> \qquad\qquad \vdots                                                                                                  &  \> \qquad \vdots \\
    \begin{array}{l}
      \mu_{\lambda_n} \cdot \omega^{n}\,\oplus\,\mu_{\lambda_n^{-1}} \cdot \omega^n,
    \end{array}                                                              & \text{ if } \tau \geq t + (n-1), \>\>\> \text{and } b = 2n-1 \text { is odd},  \\    
    \qquad \qquad           \text{or}         &  \\
     \begin{array}{l}
       \mu_{\lambda_n} \cdot \omega^{n+1}\,\oplus\,\mu_{\lambda_n^{-1}} \cdot \omega^n,  \\     
      \mathrm{ind}(\omega_2^{b+1+n(p-1) }),                                                             
    \end{array} & 
    \begin{array}{lc}
       \text{if } \tau = t + (n-1)   \\
       \text{if } \tau > t+(n-1),
   \end{array}
   \text{ and } b = 2n \text{ is even},     \\      
\end{cases} 
\end{eqnarray*}
where the $\lambda_i$ for all $1 \leq i \leq n$  are constants given by 
\begin{eqnarray*}
    \lambda_i & = &  \overline{*_i \cdot \dfrac{c}{p^{i-1}} }
\end{eqnarray*}
for some `fudge factors' $*_i$, except if $b = 2n-1$ is odd and $i = n$, in which case $\lambda_n$ satisfies
\begin{eqnarray*}
  \lambda_n +  \frac{1}{\lambda_n} & = &  \overline{*_n \cdot \dfrac{c}{p^{n-1}} }.
\end{eqnarray*}
\end{conj}

\noindent Thus, on the inertia subgroup $I_{\Q_p}$,  Conjecture~\ref{zigzag} predicts the reduction is given by the  picture:
\vskip 0.5 cm

\begin{tikzpicture}[xscale = 1.8, auto=center][extra thick]
\draw [latex-latex] (0,0) -- (8.0,0);
\foreach \x in {1.2,2.9,4.9}
\draw[shift={(\x,0)},color=black] (0pt,3.5pt) -- (0pt,-3.5pt);
\node at(.5,-0.4) {\small{$\ind(\omega_2^{b+1})$}};
\node at (1.2,0.5) {$\omega^b \oplus \omega$};
\node at (2.0,-0.4){\small{$\ind(\omega_2^{b+p})$}};
\node at (2.9,0.5) {$\omega^{b-1} \oplus \omega^2$};
\node at (3.7,0.5) {$\cdots$};
\node at (3.7,-0.4) {$\cdots$};
\node at (4.9,0.6){\tiny{$\begin{cases} \omega^n \oplus \omega^n,  & b=2n-1  \\  \omega^{n+1} \oplus \omega^n, & b=2n  \end{cases} $}};
\node at (7.0,0.4) {\small{$\omega^n \oplus \omega^n, \>\> b=2n-1$}};
\node at (6.8,-0.4) {\small{$\ind(\omega_2^{b+1+n(p-1)}), \>  b=2n $}};
\node at (1.2,-0.4) {\small{$t$}};
\node at (2.9, -0.4){\small{$t+ 1$}};
\node at (4.9,-0.4){\small{$t+(n-1)$}};
\node at (8.2, 0) {$\tau$};
\end{tikzpicture}

\vskip 0.5 cm

One might further refine the conjecture and give the shape of the explicit parameter $c$. 
One possibility for $c$ that works for small
half-integral slopes is as follows:
\begin{eqnarray}
    \label{c}
    c & = &  \frac{a_p^2 - \binom{r-v_-}{v_+} \binom{r - v_+}{v_-} p^b}{p a_p},
\end{eqnarray}
where $v_-$ and $v_+$ are the largest, respectively smallest, integers not equal to $v  \in \frac{1}{2} 
{\mathbb Z}$ such that $v$ lies in the 
interval $(v_-, v_+)$,
Moreover, we expect that the fudge factors $*_i$ in the conjecture are a simple rational expression in $r$, though we 
have not been able to guess a precise formula for it in general. However,
we expect that $*_1 =  \dfrac{b}{b-r}$, and we will soon give a formula for $*_2$ when $v = \frac{3}{2}$.

\subsection{Evidence}
As mentioned earlier, the conjecture is true for $v=\frac{1}{2}$ \cite{BG13} and $v = 1$ \cite{BGR18} for all weights $r > b$,
including the refinement involving the shape of $c$ above and the constants $\lambda_i$. 
As further evidence towards the refined form of the conjecture, we present  \cite[Theorem 1.1]{GR19}:

\begin{theorem} 
  \label{maintheorem}
  If $v = \frac{3}{2}$, then the zig-zag conjecture is true.  More precisely, if $p \geq 5$, $v = \frac{3}{2}$, $r > b=3$,  and 
\begin{eqnarray*}
    c & = &  \frac{a_p^2 - (r-2) \binom{r-1}{2} p^3}{p a_p}, 
  \end{eqnarray*}
  and we set
  \begin{eqnarray*}
  \tau & = & v \left( c \right), \\
      t  & = & v (b-r),
  \end{eqnarray*}
then  the reduction $\bar{V}_{k,a_p}$ enjoys the following {\it tetrachotomy}:
\begin{eqnarray*}
\bar{V}_{k, a_p} & \sim & 
  \begin{cases}
    \mathrm{ind}(\omega_2^{b+1}),  & \text{if } \tau < t \\
    \mu_{\lambda_1} \cdot \omega^b\,\oplus\,\mu_{\lambda_1^{-1}} \cdot \omega, & \text{if } \tau = t \\ 
    \mathrm{ind}(\omega_2^{b+p}), & \text{if } t <  \tau < t+1 \\       
    \mu_{\lambda_2} \cdot \omega^{b-1}\,\oplus\,\mu_{\lambda_2^{-1}} \cdot \omega^2, & \text{if } \tau \geq t+1, 
    \end{cases}
\end{eqnarray*}
where the $\lambda_i$ are constants given by
\begin{eqnarray*}
  \lambda_1 &  = & \overline{\dfrac{b}{b-r} \cdot c } \\
  \lambda_2 + \frac{1}{\lambda_2} & = &  \overline{\dfrac{b-1}{(b-1-r){(b-r)}}\cdot \dfrac{c}{p}}.
\end{eqnarray*}
\end{theorem}

Since $\bar{V}_{k,a_p}$ was completely determined in \cite{Bhattacharya-Ghate} for all other slopes $v$ in $(1,2)$ not equal to $\frac{3}{2}$,
(even for $p \geq 3$), and in \cite{BGR18} for slope $v = 1$ (for $p \geq 5)$, we obtain the following corollary.

\begin{cor}
  If $p \geq 5$, then the reduction $\bar{V}_{k, a_p}$ is known for all slopes $v = v(a_p)$ less than $2$. 
\end{cor}

\section{Local Constancy vs Zig-Zag}

Let us now explain why we have included the caveat  
`except possibly some weights $r$ that are $p$-adically close to some small weights'
in the qualitative statement
of the zig-zag conjecture. As Theorem~\ref{maintheorem} above and the preceding historical discussion shows,
there is no such caveat for the first few half-integral
slopes $\frac{1}{2}$, $1$ and $\frac{3}{2}$. More precisely, for these slopes, it suffices to take $r > b$.\footnote{We exclude 
$r = b$ in order to exclude the
degenerate case $t = \infty$ and $\tau = v(c) = \infty$ which happens for  $a_p^2 = p^b$ when the explicit 
parameter $c$ in equation \eqref{c} vanishes. In fact, we could include $r = b$ if we exclude $a_p = \pm p^{v}$, since 
both zig-zag (now $\tau$ is finite so less than $t = \infty$) and Fontaine-Edixhoven predict the same answer, 
namely $\ind(\omega_2^{b+1})$.}
However, this caveat is required for slopes at least 2. 

In order to explain this, we recall
the following local constancy result of Berger \cite[Theorem B]{Berger12},
whose proof uses the families of trianguline representations constructed in \cite{Col08}, \cite{Che13}.

\begin{theorem}(Local Constancy for Weights)
  \label{localconstancy}
  Suppose $a_p \neq 0$ and $k > 3 \cdot v + \alpha(k-1) + 1$, where $\alpha(n) = \sum_{j\geq1} \lfloor \frac{n}{p^{j-1}(p-1)} \rfloor$. 
  Then there is\footnote{The theorem asserts the existence of $m$. 
  For some bounds on its size for small weights,
  see Bhattacharya \cite{B19}.} 
 a positive integer
 $m = m(k,a_p)$ such that $$\bar{V}_{k',a_p} \sim \bar{V}_{k,a_p},$$ for all $k' > k$ with $k' \equiv k$ mod $p^{m-1}(p-1)$.
\end{theorem}

There is an interesting interplay between local constancy and the zig-zag conjecture. In most cases both are compatible. However, sometimes 
they are not. As we shall see below, in some cases zig-zag holds forcing local constancy to fail.
And in others, local constancy holds forcing zig-zag 
to fail. It is these last kinds of cases that we exclude from the statement of the zig-zag conjecture, whence the caveat in the statement of the conjecture. Let us elaborate. 

\subsection{Fontaine-Edixhoven weights} 
First consider weights in the Fontaine-Edixhoven range, namely $k \leq p+1$ or $r = b \leq p-1$. It is known that
$\bar{V}_{k, a_p} \sim \mathrm{ind}(\omega_2^{b+1})$ for such weights.
The only exceptional weight in this range satisfying the hypothesis of Theorem~\ref{localconstancy} is $r = b = 1$ for $v = \frac{1}{2}$, 
since $$k = 2v + 2 \not> 3 \cdot v + \alpha(b+1) + 1,$$ unless $v < 1$. So assume $r = 1$ and $v = \frac{1}{2}$.
By local constancy, if $r' = k'-2 \equiv 1$ mod $p^{t'}(p-1)$, for $t'$ sufficiently large, we must have $\bar{V}_{k',a_p} \sim \bar{V}_{k,a_p} \sim \mathrm{ind}(\omega_2^2)$. 
But if $t'$ is larger than $\tau' = v \left( \frac{a_p^2 - r'p}{pa_p} \right)$, then zig-zag (theorem of \cite{BG13}) also predicts 
$\bar{V}_{k',a_p} \sim \ind(\omega_2^2)$, so there is no apparent incompatibility with local constancy at $r = 1$.

However, as remarked above, the next exceptional weight, $r = 2$ for $v = 1$ does not satisfy the bound in Theorem~\ref{localconstancy}, 
since 
\begin{eqnarray}
  \label{bound}
  4 & \not> & 3\cdot1+ 0 + 1.
\end{eqnarray} 
In fact, we can  use zig-zag
to show that {\it local constancy does not hold} in this case! Indeed, if it did, then for $r' = k'-2 \equiv 2$ mod $p^{t'}(p-1)$, with $t'$ sufficiently large, we would have 
$\bar{V}_{k',a_p} \sim \bar{V}_{k,a_p} \sim \mathrm{ind}(\omega_2^3)$ is irreducible. But if we take $a_p = p$ say, then 
$\tau' = v \left( \frac{a_p^2 - \binom{r'}{2} p^2}{pa_p} \right) = t'$, so by zig-zag
(theorem of \cite{BGR18}) we must have $\bar{V}_{k',a_p}$ is reducible ($\sim \omega^2 \oplus \omega$  on inertia), a contradiction. 

Similarly, we can use the recent work \cite{GR19} to show that {\it local constancy fails} when $r = 3$ and say
$a_p = p^{\frac{3}{2}}$ (so $v = \frac{3}{2}$ and the bound on $k$ in Theorem~\ref{localconstancy} is not satisfied). 
If local constancy were to 
hold, then for $r' = k'-2$ as above we would have $\bar{V}_{k',a_p} \sim \bar{V}_{k,a_p} \sim \mathrm{ind}(\omega_2^4)$. 
But by \eqref{c} we have, $\tau' = \left( \frac{a_p^2 - \binom{r'-1}{2} (r'-2) p^3}{pa_p} \right) = t' + \frac{1}{2}$, so by 
zig-zag (Theorem~\ref{maintheorem}), it is now known that $\bar{V}_{k',a_p} \sim \mathrm{ind}(\omega_2^{3+p})$, the next irreducible 
possibility, a contradiction!

Let us record these observations formally now, since they do seem to have been noticed before.

\begin{theorem}
   Local constancy in the weight for the reductions $\bar{V}_{k,a_p}$ may fail for small weights when $a_p \neq 0$, 
   e.g., it fails for $(k, a_p) =(4,p)$ and $(5,p^{\frac{3}{2}})$.
\end{theorem}

\noindent  It has been known for some 
time (see \cite{Berger12}) that local constancy in the weight fails when $a_p = 0$, in view of the main result of \cite{BLZ}.
Also, the above discussion answers the second question below \cite[Theorem 1.1]{B19} in the negative, namely one {\it cannot always} improve the lower bound  $3v + \alpha(k-1) + 1$ in Theorem~\ref{localconstancy}, since \eqref{bound} shows this lower bound is sharp
when $k=4$ and $v = 1$. However, this bound can be improved for other small weights $k$ and slopes $v$, see \cite[Thm 1.2]{B19}.

\subsection{Breuil weights} We now consider the next range of small weights, namely  $p \leq r = k-2 = b+(p-1) \leq 2p-2$ which we refer to
as Breuil weights, since 
the reductions $\bar{V}_{k, a_p}$ 
at these weights were investigated completely in
\cite{Bre03b} for $p$ odd. Every exceptional Breuil
weight satisfies the hypothesis of Theorem~\ref{localconstancy}, since $$k = b + p + 1= 2v+p+1 > 3 v + 1 + 1$$  
(or even $3v + 2 + 1$, for $p > 3$). 
So if $r' = k'-2 \equiv b + (p-1)$ mod $p^{t'}(p-1)$ for $t'$ sufficiently large, we must have $\bar{V}_{k',a_p} \sim \bar{V}_{k,a_p}$. 
Now Breuil showed (see \cite[Theorem 5.2.1, parts 2 and 3]{Berger11}) that, for $k = p+2$, $p+3$ and $p+4$ and slopes $v = \frac{1}{2}$, $1$ and $\frac{3}{2}$, respectively, the reduction $\bar{V}_{k,a_p}$ is isomorphic to $\mathrm{ind}(\omega_2^2)$, or
$\mu_\lambda \cdot \omega^2 \oplus \mu_{\lambda^{-1}} \cdot \omega$ with $\lambda = \overline{2 \cdot \frac{a_p}{p}}$,
or $\mathrm{ind}(\omega_2^{p+3})$, respectively.
Now the parameter $\tau' = v(c)$, with $c = c(r', a_p)$ as in \eqref{c}, is easily checked in these cases to be minimal, namely $\tau' = -\frac{1}{2}$, $0$ and $\frac{1}{2}$, 
respectively, whereas the other parameter in zig-zag namely $v(b-r') = 0$ vanishes in all three cases. An easy check using zig-zag (a theorem 
in all three cases, including slope $\frac{3}{2}$, by Theorem~\ref{maintheorem}), shows that $\bar{V}_{k',a_p}$ is exactly isomorphic to 
the above three representations, respectively (even the formula for $\lambda$ matches well).
In other words, local constancy at the Breuil weights is compatible with zig-zag for slopes $v = \frac{1}{2}, 1$ and $\frac{3}{2}$.

However, things begin to go wrong at the next exceptional Breuil weight $k = p +5$, where $b = 4$ and $v = 2$. Breuil proves
$\bar{V}_{k,a_p}  \sim \mathrm{ind}(\omega_2^{p+4})$ (\cite[Theorem 5.2.1, part 2]{Berger11}) and so is irreducible. However, with notation as above, one 
checks $\tau' = 1$ which is one more than 
$v(b-r')=0$, so zig-zag says $\bar{V}_{k',a_p} \sim \omega^3 \oplus \omega^2$ (on inertia), which is reducible! In fact 
the first `counterexamples' found by Rozensztajn to the refined form of zig-zag for $v = 2$ and $r > b$ were of this kind - they arise as just explained 
from Breuil's results and local constancy. Subsequently, the author and Rai found other Breuil weights (for $v > 2$) 
for which zig-zag doesn't hold, so when $v \geq 2$, we 
need to exclude such small weights $r = b+p-1$ (and large weights that are $p$-adically close to them) from 
the zig-zag conjecture. 
This explains the caveat in the statement of zig-zag, at least for Breuil weights.

\subsection{Higher small weights} 
Let us now turn to general small weights $r = b + m (p-1)$, for an integer $m > 0$. When $m \geq 2$, we have $r \geq 2p-1$ 
and we call such weights super-Breuil weights. As far as we are aware, the reductions $\bar{V}_{k,a_p}$ have not been determined for super-Breuil
weights,\footnote{for all slopes, though the reductions have been determined in \cite{B19} for $m=2$ and $3$ and some small slopes.} 
except in the smallest case $r = 2p-1$, where $b = 1$ and $m = 2$ (again by Breuil, though the
answer is only stated in \cite[Theorem 5.2.1, part 4]{Berger11}). Local constancy  (Theorem~\ref{localconstancy}) holds 
for this weight since  $$k = 2p+1 > 3 \cdot \frac{1}{2} + 2 + 1,$$ for all $p$, so for $r' = k'-2 \equiv 2p-1$ mod $p^{t'}(p-1)$, for $t'$ sufficiently
large, we obtain
\begin{eqnarray}
  \label{k=2p+1}
  \bar{V}_{k',a_p} \sim \bar{V}_{k,a_p} 
                                                   \sim \begin{cases}
                                                      \mathrm{ind}(\omega_2^2) & \text{if } v(a_p^2 + p) < \frac{3}{2}, \\
                                                      \mu_\lambda \cdot  \omega \oplus \mu_{{\lambda}^{-1}} \cdot \omega        & \text{if } v(a_p^2 + p) \geq \frac{3}{2},
                                                    \end{cases}                                                   
\end{eqnarray}
for an explicit $\lambda$. On the other hand $v(b-r') = 0$, and by \eqref{c} we have
\begin{eqnarray*}
  \tau' = v \left( \frac{a_p^2 - (2p-1)p}{pa_p} \right) \text{ which is }\begin{cases}  
                                                                                     < 0 & \text{if } v(a_p^2 + p) < \frac{3}{2}, \\
                                                                                    \geq 0    & \text{if } v(a_p^2 + p) \geq \frac{3}{2},
                                                                                \end{cases}
\end{eqnarray*}                                                                                   
so zig-zag (theorem of \cite{BG13}) predicts exactly the same answers as in \eqref{k=2p+1} (with the same $\lambda$).  
Thus local constancy at the first super-Breuil weight $r = 2p-1$ is compatible with zig-zag for $v = \frac{1}{2}$.

Since as far as we are aware, there are no general theorems yet giving $\bar{V}_{k,a_p}$ for other super-Breuil weights (for all slopes), 
it is not possible to directly compare local constancy at these points with zig-zag. However, recently the author and Rai found some
further numerical `counterexamples'  to the refined version of zig-zag at super-Breuil points, showing
that for slopes $v \geq 2$, such weights (and large weights $p$-adically close to them) need 
also to be excluded. 

\subsection{A bound} The question now arises as to which small weights one needs to exclude from the zig-zag pattern. To this end
we recall that the main result of \cite{BLZ} specifies
the reduction $\bar{V}_{k, a_p}$ {\it for slopes that are large compared to the weight}. This result immediately gives
us an explicit bound 
on which small weights (and large weights that are $p$-adically close to them) one must exclude from
zig-zag. We describe this now. 

The main theorem of \cite{BLZ} says that if $v >  \lfloor \frac{k-2}{p-1} \rfloor$,
then $\bar{V}_{k, a_p} \sim \mathrm{ind}(\omega_2^{k-1})$. For $r = b + m(p-1)$, where we take $b \neq p-1$ for simplicity, this 
translates to saying that if $v > m$, then $\bar{V}_{k, a_p} \sim \ind(\omega_2^{b+1+m(p-1)})$.  Let us compare this result with
the refined form of 
zig-zag which uses the explicit form of $c$ in 
\eqref{c}. Take $v = m+\frac{1}{2}$. So $v_- = m$ and $v_+ = m+1$, and $p$ divides the binomial coefficient
${r-(m+1) \choose m}$, so $\tau = (2m+1) - (m+ \frac{3}{2}) = m - \frac{1}{2}$ and we obtain that $m-1 < \tau < m$. 
Since $t = 0$ (we are taking $m$ small, so prime to $p$), zig-zag also predicts  $\bar{V}_{k, a_p} \sim \ind(\omega_2^{b+1+m(p-1)})$, so there is no inconsistency between \cite{BLZ} and zig-zag for $r = b + m(p-1)$ and $v = m + \frac{1}{2}$.

However, if we take $v = m + 1$, 
then $v_- = m$, $v_+ = m+2$, and again $p$  divides ${r-(m+2) \choose m}$, so 
$\tau = (2m+2) - (m+2) = m = t+m$, since $t =0$, so $\tau$ is integral and zig-zag predicts a reducible answer. Thus we see 
that the refined form of zig-zag is not compatible with the main result of \cite{BLZ} 
if $v \geq  \lfloor \frac{k-2}{p-1} \rfloor +  1$. 
For a fixed $v$, this says that zig-zag may not hold for small weights $r = b+m(p-1)$ with $0 < m \leq  v - 1$ (or 
possibly only for $0 < m < v-1$, for $b = p-1$), and by local constancy, for large weights which are 
$p$-adically close to these small weights.


Thus the caveat 
`except possibly some weights $r$ that are $p$-adically close to some small weights'
in the qualitative statement of zig-zag
should be taken to mean that zig-zag should hold for all exceptional weights $k$, for all $k -2 > b$, except 
possibly for some $k$ lying
in small $p$-adic disks around Breuil and super-Breuil weights of the form $r = b + m(p-1)$ with $0 < m \leq v-1$. 
Since we are excluding the case $m = 0$, this is only a caveat for slopes $v \geq 2$!

\subsection{Lag} A final remark. Numerical computations show that for very special values of $a_p$, 
the reduction $\bar{V}_{k,a_p}$ lags behind what is predicted by zig-zag {\it even if} $r = b + m(p-1)$ and
$m > v-1$. This lag occurs, for instance, when $v = 2$ and $a_p = p^2$, for $m = 2$, $3$ etc. 
Since we have not understood 
this phenomenon yet, we have also allowed the possibility of dropping 
(hopefully) finitely many values of $a_p$ of a particular slope $v$ from the qualitative statement of the zig-zag conjecture.
%

\section{Proof of Theorem~\ref{maintheorem} and beyond} The proof  of Theorem~\ref{maintheorem} uses the compatibility
between the $p$-adic and mod
$p$ Local Langlands Correspondences first introduced in \cite{Bre03b}, with respect to the process of reduction \cite{B10}. This compatibility allows one to 
reduce the reduction problem to a representation theoretic one, namely, to computing
the reduction of a $\mathrm{GL}_2(\Q_p)$-stable lattice in a 
certain unitary $\mathrm{GL}_2(\Q_p)$-Banach space. We recall the key ingredients of the argument here
for slope $v = \frac{3}{2}$ and discuss how the details of the argument should generalize 
to larger half-integral slopes $v$. In doing this, we will see how the zig-zag conjecture acquired its name. 
\subsection{Hecke operator}
  \label{Hecke}

Let $G = \mathrm{GL}_2(\Q_p)$, $K = \mathrm{GL}_2(\Z_p)$ be the standard
maximal compact subgroup of $G$ and $Z = \Q_p^\times$ be the center of $G$.
Let $R$ be a $\Z_p$-algebra and let $V = \Sym^r R^2\otimes D^s$ be the
usual symmetric power representation of $KZ$ twisted by a power of the
determinant character $D$, modeled on homogeneous polynomials of degree
$r$ in the variables $X$ and $Y$ over $R$. We denote compact induction by $\mathrm{ind}_{KZ}^{G}$.
Thus $\mathrm{ind}_{KZ}^{G} V$ consists of functions
$f : G \rightarrow V$ such that $f(hg) = h \cdot f(g)$, for all $h \in KZ$
and $g \in G$, and $f$ is compactly supported mod $KZ$.   For $g \in G$, $v \in V$, let
$[g,v] \in \mathrm{ind}_{KZ}^{G} V$ be the function with support in
${KZ}g^{-1}$ given by 
  $$g' \mapsto
     \begin{cases}
         g'g \cdot v,  \ & \text{ if } g' \in {KZ}g^{-1} \\
         0,                  & \text{ otherwise.}
      \end{cases}$$ 
Any function in $\mathrm{ind}_{KZ}^G V$ is a finite linear combination of functions of the form $[g,v]$, for $g\in G$ and $v\in V$.  
The Hecke operator $T$ is defined by its action on these elementary functions via

\begin{equation}\label{T} T([g,v(X,Y)])=\underset{\lambda\in\F_p}\sum\left[g\left(\begin{smallmatrix} p & [\lambda]\\
                                                 0 & 1
                                                \end{smallmatrix}\right),\:v\left(X, -[\lambda]X+pY\right)\right]+\left[g\left(\begin{smallmatrix} 1 & 0\\
                                                                                                                                                   0 & p
                                                \end{smallmatrix}\right),\:v(pX,Y)\right],\end{equation}
 where $[\lambda]$ denotes the Teichm\"uller representative of $\lambda\in\F_p$.        

\subsection{The mod $p$ Local Langlands Correspondence}
\label{subsectionLLC}

For $0 \leq r \leq p-1$, $\lambda \in \bar{\F}_p$ and $\eta : \Q_p^\times
\rightarrow \bar\F_p^\times$ a smooth character, let
\begin{eqnarray*}
  \pi(r, \lambda, \eta) & := & \frac{\mathrm{ind}_{KZ}^{G} \:\Sym^r\bar\F_p^2}{T-\lambda} \otimes (\eta\circ \mathrm{det})
\end{eqnarray*}
be the smooth admissible representation of ${G}$, 
known to be irreducible unless $(r,\lambda)=(0,\pm 1)$ or $(p-1,\pm 1)$, by the classification of irreducible  representations  of $G$ in characteristic $p$ in \cite{BL94, BL95, Bre03a}. For $(r,\lambda) = (0, \pm1)$, we have exact sequences
$$  0  \rightarrow \mathrm{St} \rightarrow \pi(0,1,1) \rightarrow  \bar\F_p \rightarrow 0, $$
$$0  \rightarrow \mathrm{St}  \otimes \mu_{-1}  \rightarrow \pi(0,-1,1) \rightarrow  \bar\F_p(\mu_{-1}) \rightarrow 0,$$
where the last surjective maps are induced by `total sum', respectively `alternating sum', of the values of a 
compactly supported function on the tree, and the (irreducible) Steinberg representation $\mathrm{St}$ can be taken 
to be defined by the first sequence. 
Similar exact sequences exist for $(r, \lambda) = (p-1, \pm 1)$. 
With this notation, Breuil's semisimple mod $p$ Local Langlands Correspondence (mod $p$ LLC) is given by  
\cite[Def. 1.1]{Bre03b}: 
\begin{itemize} 
  \item  $\lambda = 0$: \quad
             $\mathrm{ind}(\omega_2^{r+1}) \otimes \eta \:\:\overset{LL}\longmapsto\:\: \pi(r,0,\eta)$,
  \item $\lambda \neq 0$: \quad
             $\left( \mu_\lambda\cdot \omega^{r+1}  \oplus \mu_{\lambda^{-1}} \right) \otimes \eta 
                   \:\:\overset{LL}\longmapsto\:\:  \pi(r, \lambda, \eta)^{ss} \oplus  \pi([p-3-r], \lambda^{-1}, \eta \omega^{r+1})^{ss}$, 
\end{itemize}
where $\{0,1, \ldots, p-2 \} \ni [p-3-r] \equiv p-3-r \mod (p-1)$.


\subsection{Standard lattice}
Recall that $G =  \mathrm{GL}_2(\Q_p)$. Let $B(V_{k,a_p})$ be the unitary $G$-Banach space associated
to $V_{k,a_p}$ by the $p$-adic Local Langlands Correspondence. 
The reduction $\overline{B(V_{k,a_p})}^{ss}$ of a lattice 
in this space coincides with the image of $\bar{V}_{k,a_p}$
under the (semisimple) mod $p$ LLC defined above.
Since the mod $p$ LLC is by definition injective, it suffices to compute the reduction 
$\overline{B(V_{k,a_p})}^{ss}$.

Recall $K = \mathrm{GL}_2(\Z_p)$ and $Z= \Q_p^\times$. Let $X =
KZ \backslash G$ be the (vertices of the) Bruhat-Tits tree attached to
$G$.  The module $\mathrm{Sym}^r \bar\Q_p^2$, for $r = k-2$, carries a
natural action of $KZ$, and the projection

$$
KZ \backslash ( G \times   {\mathrm {Sym}}^{k-2} \bar\Q_p^2 ) \rightarrow KZ \backslash G = X
$$
defines a  local system on $X$. The space 
\begin{eqnarray*}  
   {\mathrm{ind}}_{KZ}^{G} \> {\mathrm {Sym}}^{k-2} \bar\Q_p^2
\end{eqnarray*}
consisting of all sections $f : G \rightarrow \mathrm{Sym}^{k-2}
\bar\Q_p^2$ of this local system which are compactly supported mod $KZ$ is a representation for $G$,
equipped with a $G$-equivariant Hecke operator $T$. 
Let $\Pi_{k,a_p}$ be the locally algebraic representation of
$G$ defined by taking the cokernel of $T-a_p$ acting on the above space of sections.\footnote{$\Pi_{k,a_p}$ is the tensor product of a smooth representation and
an algebraic representation of $G$ and is generically irreducible; in a more general setting, D. Prasad \cite{Pr01} has shown that
irreducible locally algebraic representations always have this form.} 
Let $\Theta_{k,a_p}$ be the image of the integral sections
${\mathrm{ind}}_{KZ}^{G} \> {\mathrm {Sym}}^{k-2} \bar\Z_p^2$ in
$\Pi_{k,a_p}$. 
Then $B(V_{k,a_p})$ is the completion $\hat{\Pi}_{k,a_p}$ of
$\Pi_{k,a_p}$ with respect to the lattice $\Theta_{k,a_p}$.
The completion $\hat{\Theta}_{k,a_p}$, and sometimes by
abuse of notation $\Theta_{k,a_p}$ itself, is called the standard lattice
in $B(V_{k,a_p})$.  We have $\overline{B(V_{k,a_p})}^{ss} \cong
\bar{\hat{\Theta}}_{k,a_p}^{ss} \cong \bar\Theta_{k,a_p}^{ss}$.
Thus, to compute $\bar{V}^{}_{k,a_p}$, it suffices to compute the reduction 
$\bar\Theta_{k,a_p}^{ss}$ of the standard lattice $\Theta_{k,a_p}$.

\subsection{A filtration}
We introduce some further notation. 
Let $V_r = \Sym^r \, \bar\F_p^2$ (this is just $V$ from \S \ref{Hecke} with $R = \bar\F_p$ and $s = 0$) 
denote the $(r+1)$-dimensional $\bar{\mathbb{F}}_p$-vector space of 
homogeneous polynomials $P(X,Y)$ in two variables $X$ and $Y$ of degree $r$ over
$\bar{\mathbb{F}}_p$. 
The group $\Gamma=\mathrm{GL}_2(\mathbb{F}_p)$
acts on $V_r$ by the formula 
$\left(\begin{smallmatrix} a & b\\
c & d
\end{smallmatrix}\right)\cdot P(X,Y)=P(aX+cY, bX+dY)$, and $KZ$ acts on $V_r$ via projection
to $\Gamma$, with $\left( \begin{smallmatrix} p & 0 \\ 0 & p \end{smallmatrix} \right) \in
Z$ acting trivially. By definition of the lattice $\Theta_{k,a_p}$, there is a surjection 
${\mathrm{ind}}_{KZ}^{G} \> {\mathrm {Sym}}^{k-2} \bar\Z_p^2 \twoheadrightarrow \Theta_{k,a_p}$,
which induces a surjective map
\begin{eqnarray}
   \label{bartheta}
   \mathrm{ind}_{KZ}^G V_r \twoheadrightarrow\bar{\Theta}_{k,a_p}.
\end{eqnarray}
Thus, in order to compute $\bar\Theta_{k,a_p}$, it suffices to investigate the kernel of the map \eqref{bartheta}.
Let $\theta(X,Y)=X^pY-XY^p$.
The action of $\Gamma$ on $\theta$ is via the determinant $D : \Gamma \rightarrow \mathbb{F}_p^\times$.
For each $i \geq 0$,  consider the $\Gamma$- hence $KZ$-submodule $V_r^{(i)} = \{P(X,Y) \in V_r :  \theta^i  | P \}$
of $V_r$ consisting of all polynomials divisible by $i$-copies of the $\theta$-polynomial.
Clearly $V_r^{(0)} = V_r$ and $V_r^{(1)} = V_r^*$ is the largest singular submodule of $V_r$ and
$$
V_r^{(i)} \sim \begin{cases}
                                          0, &\text{if } r< i(p+1)\\
                                          V_{r-i(p+1)}\otimes D^i, &\text{if } r\geq i(p+1),
                                         \end{cases}
$$
for $i \geq 0$. The submodules $V_r^{(i)}$ are important in computing the kernel of \eqref{bartheta} because of the following useful fact
\cite[Rem. 4.4]{BG09}: if the slope $v = v(a_p)< i$ and $r \geq i(p+1)$, then $\mathrm{ind}_{KZ}^G V_r^{(i)}$ lies in
the kernel of \eqref{bartheta}, i.e.,  the surjection \eqref{bartheta}
factors through the 
map $$\mathrm{ind}_{KZ}^G \dfrac{V_r}{V_r^{(i)}} \twoheadrightarrow\bar{\Theta}_{k,a_p}.$$
In the setting of the zig-zag conjecture, the smallest $i$ we can choose is as follows.
Recall $b = 2v$ and $b = 2n-1$ is odd or $b = 2n$ is even, with $n \geq 1$.
If the former case, $v = n - \frac{1}{2}$, so we can take $i = n$, whereas in the latter case,
we have $v = n$, and so we can take $i = n+1$.

In general, the submodules $V_r^{(i)}$ for $i \geq 0$ define a filtration on $V_r$
$$V_r \supset V_r^{(1)} \supset \cdots \supset V_r^{(i)} \supset V_r^{(i+1)} \supset \cdots, $$ 
with each subquotient  $V_r^{(i)} / V_r^{(i+1)}$ containing two explicit Jordan-H\"older (JH) factors
$J_{2i}$ and $J_{2i+1}$ fitting into the exact sequence
\begin{eqnarray}
  \label{ps}
  0 \rightarrow J_{2i} \rightarrow V_r^{(i)} / V_r^{(i+1)} \rightarrow J_{2i+1} \rightarrow 0.
\end{eqnarray}
In the setting of the zig-zag conjecture, these JH factors can be described explicitly as follows: $$J_{2i} = V_{b-2i} \otimes D^{i} \quad \text{ and } \quad   J_{2i+1} = V_{p-1-b+2i} \otimes D^{b-i},$$ for 
$0 \leq i \leq n-1$ if $b=2n-1$ is odd, and $0\leq i \leq n$ if $b = 2n$ is even.
We remark that when $b = 2n$ is even, the last JH-factor $J_{2n+1} = V_{p-1} \otimes D^n$ is projective and the exact
sequence \eqref{ps} for $i = n$ splits, so we can flip the places of $J_{2n+1}$ and $J_{2n}$ in
the exact sequence \eqref{ps}
above.
Thus the associated graded of the relevant part of the filtration above fits into the following picture:
\vskip 0.2 cm
\quad
\begin{xymatrix}{
         J_0  \ar[d]  &   J_2 \ar[d]  &   J_ 4 \ar[d]   &    &   & J_{2n-2} \ar[d]  & J_{2n} \text{ or } J_{2n+1} \ar[d]  \\ 
   \dfrac{V_r}{V_r^{(1)}} \ar[d]  & \dfrac{V_r^{(1)}}{V_r^{(2)}} \ar[d] & \dfrac{V_r^{(2)}}{V_r^{(3)}}  \ar[d] & \cdots & \cdots &  \dfrac{V_r^{(n-1)}}{V_r^{(n)}}  \ar[d]  
   &   \dfrac{V_r^{(n)}}{V_r^{(n+1)}}, \ar[d]     \\
         J_1  \ar@{<-->}[uur]^{d} &   J_3 \ar@{<-->}[uur]^d  &   J_ 5 \ar@{<-->}[uur]^d   &     &  \ar@{<-->}[uur]^d & J_{2n-1} \ar@{<-->}[uur]^d 
         \ar@{-->}@(dl,dr)_d & J_{2n+1} \text{ or } J_{2n}   \\ }  
\end{xymatrix}
\vskip 0.2 cm
\begin{center}
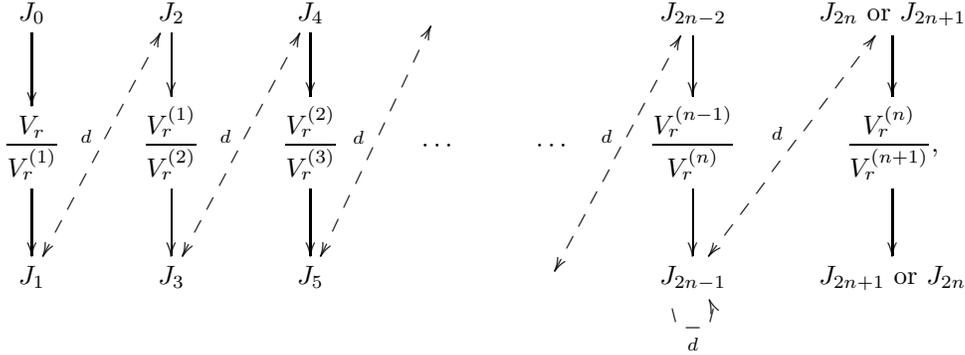
 \!\!\! A zig-zag pattern among the JH factors
\end{center} 
\vskip 0.2 cm

\noindent where the picture is taken to end at the column for ${V_r^{(n-1)}}/{V_r^{(n)}}$ for $b = 2n-1$ odd, 
and at the column for ${V_r^{(n)}}/{V_r^{(n+1)}}$ for $b = 2n$ even.  

\subsection{Proof} Recall that the mod $p$ LLC
essentially says that irreducible Galois representations 
$\bar{V}_{k,a_p}$ correspond to supersingular representations of the form $\frac{\ind_{KZ}^G J}{T}$, for some irreducible 
$\Gamma$-modules $J$, whereas reducible $\bar{V}_{k,a_p}$ correspond (generically) to a sum of two principal series representations of
the form $\frac{\ind_{KZ}^G J}{T-\lambda}$ and $\frac{\ind_{KZ}^G J'}{T-\lambda^{-1}}$, for some `dual' 
irreducible (possibly equal) $\Gamma$-modules $J$, $J'$,
the sum of whose dimensions is $p-1$ mod $(p-1)$, and some $\lambda \in \bar\F_p^\times$.  
Thus, in the picture above, exactly one, or possibly exactly two, of the above JH factors contribute to $\bar\Theta_{k,a_p}$. 
Moreover:
\begin{itemize}
\item the sum of the dimensions of $J_{2i}$ and $J_{2i+1}$ in each vertical column above is $p+1$, and by the mod $p$ LLC
each of  $J_{2i}$ and $J_{2i+1}$, if it occurs as the sole contributing factor to $\bar\Theta_{k,a_p}$, gives 
the same {\it irreducible} Galois representation for $\bar{V}_{k,a_p}$.
\item the sum of the dimensions in each adjacent diagonal pair $(J_{2i+1}, J_{2i+2})$ is $p - 1$ 
and a potential `duality' occurs (indicated
by a `$d$' on the diagonal dotted arrows), since these two JH factors may contribute together to $\bar\Theta_{k,a_p}$, giving a {\it reducible} Galois representation $\bar{V}_{k,a_p}$. 
\item when $b = 2n-1$ is odd, the last JH factor $J_{2n-1} = V_{p-2} \otimes D^n$ is potentially `self-dual' (indicated by a `$d$' on the looped dotted arrow), since twice its dimension is $0$ mod $(p - 1)$. 
\item when $b = 2n$ is even, the JH factor $J_{2n-1} = V_{p-3} \otimes D^{n+1}$ has two options with which to set up a `duality',
namely with $J_{2n}=V_0 \otimes D^n$ or with $J_{2n+1} = V_{p-1} \otimes D^n$, and in practice it does indeed `break rank' 
and pair with 
$J_{2n+1}$ sometimes (which is why we have allowed either JH factor in the last column to be the target of the last dotted
diagonal arrow `$d$'). 
\end{itemize}

We are now ready to make the key observation of this section. 
We claim that as $\tau$ varies through the rational line, the JH factors that contribute to $\bar\Theta_{k,a_p}$ actually 
occur in a  {\it zig-zag fashion} in the picture above, first down starting from $J_1$ (it is known that $J_0$ never contributes), then diagonally up and across to $J_2$ (where both $J_1$ and $J_2$ contribute together),
then only $J_2$ contributes, then down again to $J_3$ where it only contributes, then diagonally across and up to $J_4$ (where both $J_3$ and $J_4$ contribute), then down to $J_4$ where it only contributes, and so on and so forth,
as displayed by the zig-zag pattern in the diagram above. Moreover, we claim that the
diagonal jumps occur exactly when $\tau$ takes on {\it integer} values between $t$ and $t+(n-1)$ inclusive. This claim is 
rather remarkable considering that {\it a priori} there is 
no reason to expect that there should be any patterns in the way $\bar\Theta_{k,a_p}$ `selects' JH factors in the 
picture above. It also explains
why the {\it zig-zag conjecture} (Conjecture~\ref{zigzag}) has been christened as such. Finally, it explains
why the conjecture predicts that the reduction $\bar{V}_{k,a_p}$ alternates between irreducible and reducible possibilities, with the reducible possibilities occurring exactly at the aforementioned integer points (and for $\tau \geq n-1$, when $b = 2n-1$ is odd).

All of this is best summarized with another picture. Let $F_i$ for $i \geq 0$ be the sub-quotient of $\bar\Theta_{k,a_p}$ occurring
as the image of $\ind_{KZ}^G J_i$ for $i \geq 0$. Then we expect that all $F_i = 0$  vanish in 
$\bar\Theta_{k,a_p}$, {\it except} for the following
$F_i$, occurring exactly when $\tau$ is in the following regions:
\vskip 0.2 cm
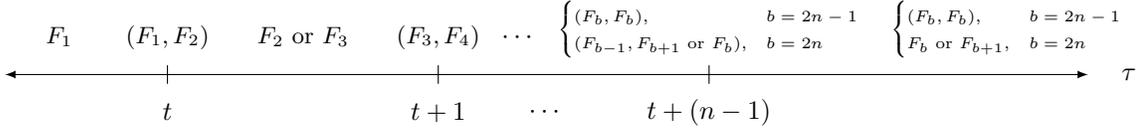
\begin{tikzpicture}[xscale = 1.8, auto=center][extra thick]
\draw [latex-latex] (0,0) -- (8,0);
\foreach \x in {1.2,3.2,5.2}
\draw[shift={(\x,0)},color=black] (0pt,3.5pt) -- (0pt,-3.5pt);
\node at(.4,0.5) {\small{$F_1$}};
\node at (1.2,0.5) {\small{$(F_1,F_2)$}};
\node at (2.2,0.5){\small{$F_2$ or $F_3$}};
\node at (3.2,0.5) {\small{$(F_3,F_4)$}};
\node at (3.8,0.5) {$\cdots$};
\node at (4.0,-0.5) {$\cdots$};
\node at (5.2,0.6){\tiny{$\begin{cases} (F_b,F_b),  & b=2n-1  \\ (F_{b-1}, F_{b+1} \text{ or } F_{b}), & b=2n  \end{cases} $}};
\node at (7.4,0.6) {\tiny{$\begin{cases} (F_b,F_b),  & b=2n-1  \\ F_{b} \text{ or } F_{b+1} , & b=2n  \end{cases} $}};
\node at (1.2,-0.5) {$t$};
\node at (3.2, -0.5){$t+ 1$};
\node at (5.2,-0.5){$t+(n-1)$};
\node at (8.3, 0) {$\tau$};
\end{tikzpicture}

%
%

In particular, when $v = \frac{3}{2}$ and $b = 3$, zig-zag predicts that the
$F_i$, for $i =1$, $2$, $3$, of $\bar\Theta_{k,a_p}$ occur in the above order.
Indeed, one might even expect the following (slightly more refined) picture holds:
\vskip 0.2 cm
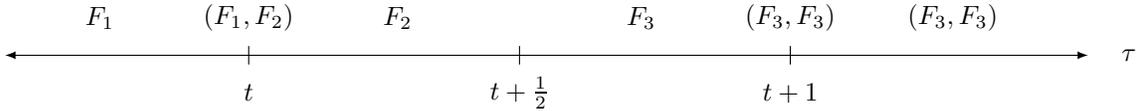
\begin{tikzpicture}[xscale = 1.8, auto=center][extra thick]
\draw [latex-latex] (0,0) -- (8,0);
\foreach \x in {1.8,3.8,5.8}
\draw[shift={(\x,0)},color=black] (0pt,3.5pt) -- (0pt,-3.5pt);
\node at(.7,0.5) {$F_1$};
\node at (1.8,0.5) {$(F_1,F_2)$};
\node at (2.9,0.5){$F_2$};
\node at (4.7,0.5) {$F_3$};
\node at (5.8,0.5) {$(F_3,F_3)$};
\node at (7,0.5){$(F_3,F_3)$};
\node at (1.8,-0.5) {$t$};
\node at (3.8, -0.5){$t+ \frac{1}{2}$};
\node at (5.8,-0.5){$t+1$};
\node at (8.3, 0) {$\tau$};
\end{tikzpicture} 
\vskip 0.1 cm
\noindent This is proved in \cite{GR19}. More precisely, using delicate computations with the Hecke operator $T$,
the following nine symmetric statements are established
in \cite[Proposition 1.2]{GR19}:
  \begin{enumerate}
  \item Around $t$:
  \begin{itemize}
    \item $\tau > t \implies F_1 = 0$
    \item $\tau = t \implies F_1 \twoheadleftarrow \dfrac{\ind J_1}{T - \lambda_1^{-1}}$ and 
                                      $F_2 \twoheadleftarrow \dfrac{\ind J_2}{T - \lambda_1}$,
                                        with $\lambda_1  = \overline{\dfrac{b}{b-r} \cdot c }$
    \item $\tau < t \implies F_2 = 0$,
 \end{itemize}
 \item Around $t + \frac{1}{2}$:
 \begin{itemize}
    \item $\tau > t + \frac{1}{2} \implies F_2 = 0$
    \item $\tau = t + \frac{1}{2} \implies F_2 \twoheadleftarrow \dfrac{\ind J_2}{T}$ and $F_3 \twoheadleftarrow \dfrac{\ind J_3}{T}$
    \item $\tau < t + \frac{1}{2} \implies F_3 = 0,$\footnote{In fact, we only prove this for $\tau \leq t$, but this suffices.}
  \end{itemize}  
  \item Around $t + 1$:
  \begin{itemize}
    \item $\tau > t + 1 \implies  F_3 \twoheadleftarrow \dfrac{\ind J_3}{T^2+1}$
    \item $\tau = t + 1 \implies  F_3 \twoheadleftarrow \dfrac{\ind J_3}{T^2-dT+1}$,
                                               with 
                                                      $d =   \overline{\dfrac{b-1}{(b-1-r)(b-r)} \cdot \dfrac{c}{p}}$.   
    \item $\tau < t + 1 \implies F_3  \twoheadleftarrow \dfrac{\ind J_3}{T}$,
  \end{itemize}  
  \end{enumerate}                                                                  
where we have written `$\ind$' for `$\ind_{KZ}^G$' for simplicity. 
Theorem~\ref{maintheorem} now follows immediately from these nine statements 
and the mod $p$ LLC, proving  zig-zag holds for all $r > b$ when the slope $v= \frac{3}{2}$.

While higher cases of the zig-zag conjecture (Conjecture~\ref{zigzag}) have yet to be attempted (it required a long 
paper just to write up all the details of the proof of the case of slope $\frac{3}{2}$), 
we expect that the selection of the JH factors that go into the proof will follow the general zig-zag pattern outlined above.

\section{Compatibility}

\subsection{Irreducibility conjecture} 
  There is a general conjecture, attributed to Buzzard, Breuil and Emerton  \cite[Conjecture 4.1.1]{BG16}, which says that 
  $\bar{V}_{k,a_p}$ is irreducible if $k$ is even and $v = v(a_p)$ is fractional, i.e., non-integral.  
  All computations of $\bar{V}_{k,a_p}$ so far support this conjecture. We remark that 
  the zig-zag conjecture is also (vacuously!) consistent with this conjecture, since the weight 
  $k$ is even exactly when the residue class $b=2v$ is even, and this happens only if $v$ is integral. 
  
\subsection{Theta operator}
  The author has suspected for some time (e.g., see the slides of his talk in the Fields symposium in 2016 in honor of 
  Bhargava) 
  that some of the reductions $\bar{V}_{k,a_p}$ in slope $v+1$ should be 
  related to the reductions in slope $v$ by twisting by $\omega$. Some evidence for this would be provided 
  if for a given form $f$ of level $N$ coprime to $p$ and slope $v$, the twisted (global) representation 
  $\bar\rho_f \otimes \omega$ is isomorphic to $\bar\rho_g$ for some form $g$ of level $N$ coprime to $p$ and slope $v' = v+1$.
  That this is indeed true (under some assumptions)
  was recently proved\footnote{That a cusp form $g$ of {\it some slope} exists follows from the seminal work of Khare-Wintenberger 
  \cite{KW09} and Kisin \cite{Kis09} on Serre's modularity conjecture. However, the forms $g$ of minimal Serre weight arising in 
  Serre's conjecture do not necessarily have slope $v+1$.}
  using the $\theta$-operator acting on overconvergent $p$-adic modular forms and the theory of
  Hida-Coleman families in \cite{GK}.  Let us compare this result with the refined version of zig-zag.
  Generically for $f$, we have $t = 0$ and $\tau$ takes its minimal value which by \eqref{c} is $2v - (v+1) = v -1$.
  So by zig-zag (cf. Conjecture~\ref{zigzag}), we have, for $n \geq 1$,
  \begin{eqnarray}
  \label{f}
     \bar\rho_f |_{I_{\Q_p}} \simeq  \begin{cases}
                                                      \ind(\omega_2^{b+1 + (n-1)(p-1)}),       & \text{if } v=n-\frac{1}{2} \\
                                                      \omega^{b-n+1} \oplus \omega^{n},   & \text{if } v=n, 
                                               \end{cases}
  \end{eqnarray} 
  depending on whether $v$ is half-integral or integral. 
  It turns out that the weight $l$ of $g$ satisfies $l \equiv k + 2 \mod (p-1)$, so $l$ is an exceptional weight for $g$ 
  since $l  \equiv 2v' +2  \mod (p-1)$ if and only if  $k \equiv 2v  + 2 \mod (p-1)$. Let $b' = b + 2 \in \{1, 2, \ldots, p-1 \}$ be the 
  residue class of $l$ mod $(p-1)$. We apply zig-zag to $g$. Generically, one of the parameters in zig-zag,
  namely $v(b' - (l-2)) = 0$ 
  vanishes,  whereas generically the other parameter $\tau' = v'-1 = v$  is minimal again.
  By Conjecture~\ref{zigzag}, we have
  \begin{eqnarray} 
  \label{g}
  \bar\rho_g |_{I_{\Q_p}} \simeq  \begin{cases}
                                                      \ind(\omega_2^{b'+1 + n(p-1)}),       & \text{if } v'  = n+\frac{1}{2} \\
                                                      \omega^{b'-n} \oplus \omega^{n+1},   & \text{if } v ' =n+1, 
                                               \end{cases}
   \end{eqnarray} 
  again depending on whether  the slope $v'$ of $g$ is half-integral or integral.                                           
  But as the reader may easily check, the expressions \eqref{f} and \eqref{g}
  for $\bar\rho_f |_{I_{\Q_p}}$ and $\bar\rho_g |_{I_{\Q_p}}$ are exactly compatible with the fact that the 
  local representation
  $\bar\rho_g |_{I_{\Q_p}}$ is isomorphic to $\bar\rho_f |_{I_{\Q_p}}$ twisted by $\omega$.

\end{document}